\documentclass [twoside,reqno,  12pt] {amsart}

\usepackage{amsfonts}
\usepackage{amssymb}
\usepackage{a4}
\usepackage{color}

\newtheorem{thm}{Theorem}[section]

\newtheorem{rem}[thm]{Remark}

\theoremstyle{definition}

\numberwithin{equation}{section}

\newcommand{\C}{\mathbb{C}}

\newcommand{\R}{\mathbb{R}}

\newcommand{\supp}{\operatorname{supp}}

\def \bfo {\begin {eqnarray*} }
\def \efo {\end {eqnarray*} }
\def \ba {\begin {eqnarray*} }
\def \ea {\end {eqnarray*} }
\def \beq {\begin {eqnarray}}
\def \eeq {\end {eqnarray}}
\def \supp {\hbox{supp }}

\def \p {\partial}

\def \bfo {\begin {eqnarray*} }
\def \efo {\end {eqnarray*} }
\def \ba {\begin {eqnarray*} }
\def \ea {\end {eqnarray*} }
\def \beq {\begin {eqnarray}}
\def \eeq {\end {eqnarray}}
\def \supp {\hbox{supp }}

\def \p {\partial}

\begin{document}

 \title[A remark on partial data inverse problems for elliptic equations]{A remark on partial data inverse problems for semilinear elliptic equations}

\author[Krupchyk]{Katya Krupchyk}

\address
        {K. Krupchyk, Department of Mathematics\\
University of California, Irvine\\
CA 92697-3875, USA }

\email{katya.krupchyk@uci.edu}

\author[Uhlmann]{Gunther Uhlmann}

\address
       {G. Uhlmann, Department of Mathematics\\
       University of Washington\\
       Seattle, WA  98195-4350\\
       USA\\
        and Institute for Advanced Study of the Hong Kong University of Science and Technology}
\email{gunther@math.washington.edu}

\maketitle

\begin{abstract}
We show that the knowledge of the Dirichlet--to--Neumann map on an arbitrary open portion of the boundary of a domain in $\R^n$, $n\ge 2$, for a class of semilinear elliptic equations, determines the nonlinearity uniquely.  
\end{abstract}

\section{Introduction and statement of results}
The work \cite{KLU18} discovered that nonlinearity can be helpful in solving inverse problems for hyperbolic PDE, see also \cite{CLOP}, \cite{LUW}, and the references given there. Similar phenomena for inverse problems for semilinear elliptic PDE have been revealed in  the recent works \cite{Feizmohammadi_Oksanen}, \cite{LLLS}.  A common feature of all of the aforementioned works is that the presence of a nonlinearity allows one to solve inverse problems for non-linear equations in cases where the corresponding inverse problem in the linear setting is open.   

The purpose of this note is to point out that the same phenomenon remains valid for partial data inverse boundary problems for a class of semilinear elliptic PDE.  To state our result, let $\Omega\subset \R^n$, $n\ge 2$, be a connected bounded open set with $C^\infty$ boundary. Following \cite{Feizmohammadi_Oksanen},  we consider a function  $V:\overline{\Omega}\times \C\to \C$ satisfying the conditions: 
\begin{itemize}
\item[(i)] the map $\C\ni z\mapsto V(\cdot,z)$ is holomorphic with values in the H\"older space  $C^\alpha(\overline{\Omega})$, for some $0<\alpha<1$, 
\item[(ii)] $V(x,0)=\p_z V(x,0)=0$, for all $x\in \overline{\Omega}$.
\end{itemize}
It follows from (i) and (ii) that $V$ can be expanded into a power series 
\begin{equation}
\label{eq_V}
V(x,z)=\sum_{k=2}^\infty V_k(x) \frac{z^k}{k!},\quad V_k(x):=\p_z^k V(x,0)\in C^\alpha(\overline{\Omega}),
\end{equation}
converging in the $C^\alpha(\overline{\Omega})$ topology.

Let us consider the Dirichlet problem for the following semilinear elliptic equation, 
\begin{equation}
\label{eq_ref_1}
\begin{cases}
-\Delta u+V(x,u)=0 \quad \text{in}\quad \Omega, \\
u=f \quad \text{on}\quad \p \Omega. 
\end{cases}
\end{equation}
It was shown in ~\cite{Feizmohammadi_Oksanen}, \cite{LLLS}, see also \cite{Gil_Tru_book}, that  there exist $r_0, r_1>0$ sufficiently small such that 
when  $f\in B_{r_0}^\alpha (\p \Omega):=\{f\in C^{2,\alpha}(\p \Omega): \|f\|_{C^{2,\alpha}(\p \Omega)}\le r_0\}$, the problem \eqref{eq_ref_1} has a unique solution $u=u_f\in B_{r_1}^\alpha (\overline{\Omega}):=\{ u\in C^{2,\alpha}(\overline{\Omega}): \|u\|_{C^{2,\alpha}(\overline{\Omega})}\le r_1\}$. Moreover, there is a constant $C>0$ depending on $r_0$, $r_1$ only such that 
\[
\|u\|_{C^{2,\alpha}(\overline{\Omega})}\le C \|f\|_{C^{2,\alpha}(\p \Omega)},
\]
for all $f\in B_{r_0}^\alpha (\p \Omega)$.

Let $\Gamma\subset \p \Omega$ be an arbitrary non-empty open subset of the boundary $\p \Omega$. Associated to the problem \eqref{eq_ref_1}, we define the partial Dirichlet--to--Neumann map $\Lambda_V^{\Gamma}f = \partial_{\nu} u|_{\Gamma}$, where $f\in B_{r_0}^\alpha (\p \Omega)$, $\supp(f)\subset \Gamma$. Here $u \in B_{r_1}^\alpha (\overline{\Omega})$ is the unique solution of \eqref{eq_ref_1} and $\nu$ is the unit outer normal to the boundary.

The main result of this note is as follows.
\begin{thm}
\label{thm_main}
Let $\Omega\subset \R^n$, $n\ge 2$, be a connected bounded open set with $C^\infty$ boundary, and let $\Gamma\subset \p \Omega$ be an arbitrary open non-empty subset of the boundary $\p \Omega$.  Let $V^{(1)}, V^{(2)}: \overline{\Omega}\times \C\to \C$ satisfy the assumptions (i) and (ii).  Assume that $\Lambda_{V^{(1)}}^{\Gamma} f=\Lambda_{V^{(2)}}^{\Gamma}f$ for all $f\in B_{r_0}^\alpha (\p \Omega)$ with $\supp(f)\subset \Gamma$. Then $V^{(1)}=V^{(2)}$ in $\overline{\Omega}\times \C$.
\end{thm}

\begin{rem}  In particular Theorem \ref{thm_main} shows that the potential $q\in C^\alpha(\overline{\Omega})$ in the semilinear Schr\"odinger equation with a power nonlinearity, 
\[
-\Delta u +q(x)u^m=0 \quad \text{in}\quad \Omega,
\] 
is uniquely determined by the partial Dirichlet--to--Neumann map for $m\ge 2$, when both the data and the measurements are confined to an arbitrary open portion $\Gamma\subset \p \Omega$. It may be interesting to note that the corresponding partial data inverse problem for the linear Schr\"odinger equation, i.e. $m=1$, is still open in dimensions $n \geq 3$. We refer to \cite{IUY_2010} and \cite{IY_2013} for the study of the partial data inverse problem in the linear and non-linear settings, respectively,  in dimension $n=2$.
\end{rem}

\begin{rem}
Theorem \ref{thm_main} is an immediate consequence of the main result of~\cite{DKSU} combined with the higher order linearization procedure introduced in \cite{Feizmohammadi_Oksanen},  \cite{LLLS}.
\end{rem}

Let us finally mention that inverse problems for nonlinear elliptic PDE have been  studied extensively, both in the semilinear setting, see  \cite{IY_2013},  \cite{IsaNach_1995},  \cite{IsaSyl_94}, \cite{Sun_2010}, as well as the 
quasilinear  one, see \cite{CNV_2019},  \cite{Kang_Nak_02}, \cite{Mun_Uhl_2018}, \cite{Sun_96}, \cite{Sun_Uhlm_97}.  In particular, the second order linearization of the nonlinear Dirichlet--to--Neumann map has already been used in the works   \cite{CNV_2019}, \cite{Kang_Nak_02}, \cite{Sun_96}, and \cite{Sun_Uhlm_97}.

\section{Proof of Theorem \ref{thm_main}}
\label{sec_1}

Let $\varepsilon=(\varepsilon_1, \dots, \varepsilon_m)\in \C^m$, $m\ge 2$,  and consider the Dirichlet problem \eqref{eq_ref_1} with 
\[
f=\sum_{k=1}^m \varepsilon_k f_k, \quad f_k\in C^{\infty}(\p \Omega), \quad \supp(f_k)\subset \Gamma, \quad k=1, \dots, m. 
\]
Then for all $|\varepsilon|$ sufficiently small, the problem \eqref{eq_ref_1} has a unique solution $u(\cdot,\varepsilon)\in B_{r_1}^\alpha (\overline{\Omega})$ which is holomorphic in a neighborhood of $\varepsilon =0$ in the $C^{2,\alpha}(\overline{\Omega})$ topology, see  \cite{Feizmohammadi_Oksanen}.  

Following \cite{Feizmohammadi_Oksanen}, we use an induction argument on $m\ge 2$ to prove that all the coefficients $V_m(x)$ of $V$ in \eqref{eq_V} can be determined from the partial Dirichlet--to--Neumann map.  The computations  below  related to the higher order linearization procedure reproduce those of \cite{Feizmohammadi_Oksanen}, \cite{LLLS}, the only minor difference being that we consider the case of partial measurements.

First let $m=2$, and following \cite{Feizmohammadi_Oksanen},  \cite{LLLS}, we shall proceed to carry out a second order linearization of the partial Dirichlet--to--Neumann map. Let $u_j=u_j(x; \varepsilon)\in B_{r_1}^\alpha (\overline{\Omega})$ be the unique small solution of the Dirichlet problem
\begin{equation}
\label{eq_1_1}
\begin{cases}
-\Delta u_j+ \sum_{k=2}^\infty V_k^{(j)} (x)\frac{u_j^{k}}{k!}=0 \quad \text{in}\quad \Omega, \\
u_j=\varepsilon_1 f_1+\varepsilon_2 f_2 \quad \text{on}\quad \p \Omega,
\end{cases}
\end{equation}
for $j=1,2$.   Differentiating \eqref{eq_1_1} with respect to $\varepsilon_l$, $l=1,2$, and using that $u_j(x,0)=0$, we get 
\begin{equation}
\label{eq_1_2}
\begin{cases}
\Delta v_j^{(l)}=0 \quad \text{in}\quad \Omega, \\
v_j^{(l)} = f_l \quad \text{on}\quad \p \Omega,
\end{cases}
\end{equation}
where $v_j^{(l)}= \p_{ \varepsilon_l} u_j|_{\varepsilon=0}$. By the uniqueness and the elliptic regularity for the Dirichlet problem \eqref{eq_1_2}, we see that  $v^{(l)}:=v_1^{(l)}=v_2^{(l)}\in C^\infty(\overline{\Omega})$, $l=1,2$. 

Applying $\p_{\varepsilon_1}\p_{\varepsilon_2}|_{\varepsilon=0}$ to \eqref{eq_1_1}, we get
\begin{equation}
\label{eq_1_3}
\begin{cases}
-\Delta (\p_{\varepsilon_1}\p_{\varepsilon_2} u_j |_{\varepsilon=0}) + V^{(j)}_2(x) \p_{\varepsilon_1} u_j|_{\varepsilon=0} \p_{\varepsilon_2} u_j|_{\varepsilon=0} =0 \quad \text{in}\quad \Omega, \\
\p_{\varepsilon_1}\p_{\varepsilon_2} u_j |_{\varepsilon=0} = 0 \quad \text{on}\quad \p \Omega,
\end{cases}
\end{equation}
since $\p_{\varepsilon_1}\p_{\varepsilon_2} (\sum_{k=3}^\infty V_k^{(j)} (x)\frac{u_j^{k}}{k!})$ is a sum of terms each of them containing positive powers of $u_j$, which vanish when $\varepsilon=0$, and the only term in $\p_{\varepsilon_1}\p_{\varepsilon_2} (V_2^{(j)} (x)\frac{u_j^{2}}{2!})$ which does not contain a positive power of $u_j$ is 
$V^{(j)}_2(x) \p_{\varepsilon_1} u_j  \p_{\varepsilon_2} u_j$. 
Setting $w_j=  \p_{\varepsilon_1} \p_{\varepsilon_2} u_j|_{\varepsilon=0}$, we get from \eqref{eq_1_3} that 
\begin{equation}
\label{eq_1_4}
\begin{cases}
-\Delta w_j+ V_2^{(j)}(x) v^{(1)}v^{(2)}=0 \quad \text{in}\quad \Omega, \\
w_j = 0 \quad \text{on}\quad \p \Omega.
\end{cases}
\end{equation}
The fact that $\Lambda_{V^{(1)}}^{\Gamma}(\varepsilon_1 f_1+\varepsilon_2 f_2)=\Lambda_{V^{(2)}}^{\Gamma}(\varepsilon_1 f_1+\varepsilon_2 f_2)$ for all small $\varepsilon_1, \varepsilon_2$ and all $f_1, f_2\in C^\infty(\p \Omega)$ with $\supp(f_1), \supp(f_2)\subset \Gamma$ implies that $\p_\nu u_1|_{\Gamma}=\p_\nu u_2|_{\Gamma}$.  Therefore, an application of $\p_{\varepsilon_1}\p_{\varepsilon_2}|_{\varepsilon=0}$ gives $\p_\nu w_1|_{\Gamma}=\p_\nu w_2|_{\Gamma}$. Multiplying \eqref{eq_1_4} by $v^{(3)} \in C^{\infty}(\overline{\Omega})$ harmonic in $\Omega$ and applying Green's formula,  we get
\[
\int_{\Omega}(V_2^{(1)}-V_2^{(2)})v^{(1)}v^{(2)}v^{(3)}dx=\int_{\p \Omega\setminus \Gamma} (\p_\nu w_1- \p_\nu w_2)v^{(3)}dS=0,
\]
provided that $\supp(v^{(3)}|_{\partial \Omega})\subset \Gamma$. Hence, we obtain that
\[
\int_{\Omega}(V_2^{(1)}-V_2^{(2)})v^{(1)}v^{(2)}v^{(3)}dx=0
\]
for any $v^{(l)}\in C^{\infty}(\overline{\Omega})$ harmonic in $\Omega$, such that $\supp(v^{(l)}|_{\partial \Omega})\subset \Gamma$, $l=1,2,3$.

Now taking  $v^{(3)}\not\equiv 0$ and using the result of \cite{DKSU} which says that the set of products of two harmonic functions in $C^\infty(\overline{\Omega})$ which vanish on the closed proper subset $\p \Omega\setminus \Gamma$ of the boundary $\p \Omega$ is dense in $L^1(\Omega)$, we conclude that
\[
(V_2^{(1)}-V_2^{(2)})v^{(3)}=0 \quad \text{in}\quad \Omega.
\]
Now $v^{(3)}$ is harmonic and therefore,  the set $(v^{(3)})^{-1}(0)$ is of measure zero, see \cite{Mityagin}. Hence $V_2^{(1)}=V_2^{(2)}$ in $\Omega$. 

Let $m\ge 3$ and assume that $V_k:=V_k^{(1)}=V_k^{(2)}$  in $\Omega$ for all $k=2,\dots, m-1$.  To show that  $V_m^{(1)}=V_m^{(2)}$, following  \cite{Feizmohammadi_Oksanen},  \cite{LLLS}, we shall perform the $m$th order linearization of the partial Dirichlet--to--Neumann map. To that end,  let $u_j=u_j(x; \varepsilon)\in B_{r_1}^\alpha (\overline{\Omega})$ be the unique small solution of the Dirichlet problem
\begin{equation}
\label{eq_1_5}
\begin{cases}
-\Delta u_j+ \sum_{k=2}^\infty V_k^{(j)} (x)\frac{u_j^{k}}{k!}=0 \quad \text{in}\quad \Omega, \\
u_j=\varepsilon_1 f_1+\dots +\varepsilon_m f_m \quad \text{on}\quad \p \Omega,
\end{cases}
\end{equation}
for $j=1,2$.  Next we would like to apply $\p_{\varepsilon_1}\dots \p_{\varepsilon_m}|_{\varepsilon=0}$ to \eqref{eq_1_5}. 
We first observe that  $\p_{\varepsilon_1}\dots \p_{\varepsilon_m} (\sum_{k=m+1}^\infty V_k^{(j)} (x)\frac{u_j^{k}}{k!})$ is a sum of terms each of them containing positive powers of $u_j$, which vanish when $\varepsilon=0$. Moreover, the only term in $\p_{\varepsilon_1}\dots \p_{\varepsilon_m} (V_m^{(j)} (x)\frac{u_j^{m}}{m!})$ which does not contain a positive power of $u_j$ is 
$V^{(j)}_m(x) \p_{\varepsilon_1} u_j \dots \p_{\varepsilon_m} u_j$.  Finally, the expression 
$\p_{\varepsilon_1}\dots \p_{\varepsilon_m} (\sum_{k=2}^{m-1} V_k  (x)\frac{u_j^{k}}{k!})|_{\varepsilon=0}$ contains only derivatives of $u_j$ of the form $\p^k_{\varepsilon_{l_1}, \dots, \varepsilon_{l_k}}u_j|_{\varepsilon=0}$ with $k=1, \dots, m-1$, $\varepsilon_{l_1}, \dots, \varepsilon_{l_k}\in \{ \varepsilon_1, \dots, \varepsilon_m\}$. We claim that  $\p^k_{\varepsilon_{l_1}, \dots, \varepsilon_{l_k}}u_1|_{\varepsilon=0}= \p^k_{\varepsilon_{l_1}, \dots, \varepsilon_{l_k}}u_2|_{\varepsilon=0}$ for  $k=1, \dots, m-1$, $\varepsilon_{l_1}, \dots, \varepsilon_{l_k}\in \{ \varepsilon_1, \dots, \varepsilon_m\}$. This follows by applying the operators $\p^k_{\varepsilon_{l_1}, \dots, \varepsilon_{l_k}}|_{\varepsilon=0}$ to \eqref{eq_1_5}, using the fact that  $V_k^{(1)}=V_k^{(2)}$, $k=2,\dots, m-1$,  and the unique solvability of the Dirichlet problem for the Laplacian. Hence, $H_m(x):=\p_{\varepsilon_1}\dots \p_{\varepsilon_m} (\sum_{k=2}^{m-1} V_k  (x)\frac{u_j^{k}}{k!})|_{\varepsilon=0}$ is independent of $j=1,2$.  

Therefore, we get 
\begin{equation}
\label{eq_1_6}
\begin{cases}
-\Delta (\p_{\varepsilon_1}\dots \p_{\varepsilon_m} u_j |_{\varepsilon=0}) + V^{(j)}_m(x) \p_{\varepsilon_1} u_j|_{\varepsilon=0} \dots \p_{\varepsilon_m} u_j|_{\varepsilon=0} =-H_m \quad \text{in}\quad \Omega, \\
\p_{\varepsilon_1}\dots \p_{\varepsilon_m} u_j |_{\varepsilon=0} = 0 \quad \text{on}\quad \p \Omega.
\end{cases}
\end{equation}
Proceeding as in the case $m=2$, we see that 
\[
\int_{\Omega}(V_m^{(1)}-V_m^{(2)})v^{(1)}\cdots v^{(m+1)}dx=0
\]
for any $v^{(l)}\in C^{\infty}(\overline{\Omega})$ harmonic in $\Omega$, such that $\supp(v^{(l)}|_{\partial \Omega})\subset \Gamma$, $l=1,\dots, m+1$. Arguing as in the case $m=2$, we complete the  proof of Theorem \ref{thm_main}.

\section*{Acknowledgements}
The research of K.K. is partially supported by the National Science Foundation (DMS 1815922). The research of G.U. is partially supported by NSF and a Si-Yuan Professorship of HKUST.

\end{document}